\documentclass[12pt]{article}
\usepackage{amsfonts}
\newcommand{\oper}[2]{\newcommand{#1}{\mathop{\mathrm{#2}}\nolimits} }
\oper{\tr}{tr} \oper{\adj}{adj} \oper{\Div}{div} \oper{\ad}{ad}
\oper{\Ad}{Ad} \oper{\End}{End} \oper{\Hom}{Hom} \oper{\Aut}{Aut}
\oper{\SO}{SO} \oper{\SP}{Sp} \oper{\SU}{SU} \oper{\GL}{GL}
\oper{\T}{T} \oper{\U}{U} \oper{\id}{I} \oper{\ext}{Ext}
\oper{\rank}{rank} \oper{\diag}{Diag}

\newtheorem{corollary}{Corollary}
\newtheorem{proposition}{Proposition}
\newtheorem{theorem}{Theorem}
\newtheorem{definition}{Definition}
\newtheorem{lemma}{Lemma}
\newtheorem{remark}{Remark}

\newcommand{\bproof}{\noindent{\it Proof: }}
\newcommand{\eproof}{\  q.~e.~d. \vspace{0.2in}}
\newcommand{\CnMarsOK}{\setlength{\unitlength}{0.75in}
\begin{picture}(6,0.55)
\put(0,0){\begin{picture}(1,0)
            \put(0,0.1){\circle{0.075}}
            \put(1,0.1){\circle{0.075}}
            \put(2,0.1){\circle{0.075}}
            \put(3,0.1){\circle{0.075}}
            \put(4,0.1){\circle{0.075}}
            \put(5,0.1){\circle{0.075}}
            \put(6,0.1){\circle{0.075}}
            \multiput(0.011,0.1)(1.01,0){5}{\line(1,0){0.925}}
            \put(5,0.11){\line(1,0){1}}
            \put(5,0.09){\line(1,0){1}}
            \put(5.8,0.11){\vector(-1,0){0.1}}
            \put(5.8,0.08){\vector(-1,0){0.1}}
           \end{picture}}
\end{picture}}

\newcommand{\CnblackMarsOK}{\setlength{\unitlength}{0.75in}
\begin{picture}(6,0.55)
\put(0,0){\begin{picture}(1,0)
            \put(0,0.1){\circle*{0.075}}
            \put(1,0.1){\circle{0.075}}
            \put(2,0.1){\circle*{0.075}}
            \put(3,0.1){\circle*{0.075}}
            \put(4,0.1){\circle*{0.075}}
            \put(5,0.1){\circle{0.075}}
            \put(6,0.1){\circle{0.075}}
            \multiput(0.011,0.1)(1.01,0){5}{\line(1,0){0.925}}
            \put(5,0.11){\line(1,0){1}}
            \put(5,0.09){\line(1,0){1}}
            \put(5.8,0.11){\vector(-1,0){0.1}}
            \put(5.8,0.08){\vector(-1,0){0.1}}
           \end{picture}}
\end{picture}}

\newcommand{\CnwhiteMarsOK}{\setlength{\unitlength}{0.75in}
\begin{picture}(6,0.55)
\put(0,0){\begin{picture}(1,0)
            \put(1,0.1){\circle{0.075}}
            \put(5,0.1){\circle{0.075}}
            \put(6,0.1){\circle{0.075}}
            \put(5,0.11){\line(1,0){1}}
            \put(5,0.09){\line(1,0){1}}
            \put(5.8,0.11){\vector(-1,0){0.1}}
            \put(5.8,0.08){\vector(-1,0){0.1}}
           \end{picture}}
\end{picture}}

\newcommand{\AnblackMarsOK}{\setlength{\unitlength}{0.75in}
\begin{picture}(6,0.55)
\put(0,0){\begin{picture}(1,0)
            \put(0,0.1){\circle*{0.075}}
            \put(0.6,0.1){\circle*{0.075}}
            \put(1.2,0.1){\circle{0.075}}
            \put(1.8,0.1){\circle{0.075}}
            \put(2.4,0.1){\circle{0.075}}
            \put(3,0.1){\circle*{0.075}}
            \put(3.6,0.1){\circle{0.075}}
            \put(4.2,0.1){\circle{0.075}}
            \put(4.8,0.1){\circle*{0.075}}
            \put(5.4,0.1){\circle{0.075}}
            \multiput(0.011,0.11)(0.6,0){9}{\line(1,0){0.525}}
           \end{picture}}
\end{picture}}

\begin{document}
\title{Geometry of compact complex homogeneous spaces \\
with vanishing first Chern class}
  \author{Gueo Grantcharov}

\date{}
\maketitle

\rm
\begin{abstract}
We prove that any compact complex homogeneous space with vanishing first Chern class, after an appropriate deformation of the complex structure, admits a homogeneous
Calabi-Yau with torsion structure, provided that it also has an invariant volume form. A description of such spaces
among the homogeneous C-spaces is given as well as many examples and a classification in the
3-dimensional case. We calculate the cohomology ring of some of the examples and show that in dimension 14 there are infinitely many simply-connected spaces admitting such a structure with the same Hodge numbers and torsional Chern classes. We provide also an example
solving the Strominger's equations in heterotic string theory.\\

\end{abstract}

\section{Introduction}
 A complex homogeneous space is a complex manifold which admits a transitive action of a
 complex Lie group of biholomorphisms. In this paper we consider the geometry of compact
 complex homogeneous spaces with vanishing first Chern class. Except for the tori such spaces
 are necessarily non-K\"ahler, so we are interested in the properties of an appropriate
 Hermitian connections with torsion.

 On Hermitian manifolds,
 there is a  one-parameter family of Hermitian connections
 canonically depending on the complex structure $J$ and the Riemannian metric $g$ \cite{G1}.
Among them is the Chern connection on the holomorphic tangent bundle. In this paper, we are interested in what
physicists call the K\"ahler-with-torsion connection (a.k.a. KT connection) \cite{Strom}. It is the unique
Hermitian connection whose torsion tensor is totally skew-symmetric when 1-forms are identified to their dual
vectors with respect to the Riemannian metric. If $T$ is the torsion tensor of a KT connection, it  is
characterized by the identity \cite{G1}
\[
g(T(A, B), C)=dF(JA, JB, JC)
\]
where $F$ is the K\"ahler form; $F(A,B)=g(JA,B)$,  and $A, B, C$ are any smooth vector fields.

 As a Hermitian
connection, the holonomy of a KT connection is contained in the unitary group $\U(n)$. If the holonomy of the KT
connection is reduced to $\SU(n)$, the Hermitian structure is said to be Calabi-Yau with torsion (a.k.a. CYT).

Such geometry in a physical context was considered first by A. Strominger \cite{Strom} and C. Hull \cite{Hu}.
Since then various examples were found and this led to a conjecture \cite{GIP} that {\it any} compact complex
manifold with vanishing first Chern class admits a Hermitian metric and connection with totally skew-symmetric
torsion and (restricted) holonomy in $\SU(n)$. Counterexamples to this conjecture appear in \cite{GF}. There are
also examples of CYT connections unstable under deformations. These two features of the CYT connections are in sharp
contrast to the well known moduli theory of Calabi-Yau (K\"ahler) metrics.

The main result of this paper is that after a homogeneous deformation of the complex structure any compact complex homogeneous space with vanishing first Chern class admits CYT structure if it admits a
volume form, invariant under some transitive Lie group of transformations. Moreover the CYT metric can be chosen to be invariant under the same group.
It is interesting to compare this result with the
counterexamples in \cite{GF} which are locally homogeneous. By a result of D.Guan \cite{G} any compact complex
homogeneous space with an invariant volume is a toric bundle over a product of a complex parallelizable space and
a compact homogeneous K\"ahler space. The proof of the main result is in Section 6 and uses \cite{G} and the torus bundle construction of the CYT structures in
\cite{GGP}. We also provide a characterization of the vanishing of the first Chern class. For
the complex homogeneous spaces admitting compact transitive group of transformations G (C-spaces), it is
expressed in terms of the Koszul form and the properties of the Tits fibration. This condition could be interpreted
as the vanishing of the Koszul form on the set of complementary roots of the Lie algebra of G. Such a
characterization allows us to find many explicit examples. This is the content of Section 2. In Sections 3 and 4 we consider the existence of a CYT structures on C-spaces and compact complex parallelizable manifolds. We show in Section 4 that on compact complex parallelizable manifolds, every invariant metric is CYT, and in Section 3 that on C-spaces one can always deform the complex structure through a family of homogeneous complex structures to obtain a structure admitting homogeneous CYT metric. In Section 3.3 we completely characterize the left-invariant complex structures on $SU(2)\times SU(2)$ which admit a compatible homogeneous CYT metric. It is an open set and includes the Calabi-Eckmann complex structure. In section 5 we provide the classification of compact complex homogeneous spaces in dimension three and determine the existence of a homogeneous CYT structure there. In particular we show that they admit a CYT structure except for some complex structures on $SU(2)\times SU(2)$ which are described in Section 3. In section 7 we consider the topology of
some of the examples, and as a corollary we obtain that in  dimension 14 there are infinitely many simply-connected nonflat CYT manifolds with same Betti and Hodge numbers and torsional Chern classes. In the
last section we provide also a solution to the Strominger's equations on a locally homogeneous nilmanifold.

\section{Compact complex homogeneous spaces with vanishing first Chern class.}

\subsection{Characterization of the vanishing of the first Chern class}
 We concentrate here on
the compact complex homogeneous spaces admitting transitive action of a (real) compact Lie group.
 Such spaces are called C-spaces and were investigated first
by Wang \cite{W}.
 In \cite{W} is proved the following  result:

 \begin{theorem} Let $G$ be a compact semi-simple Lie group
 and $H$ be
a closed connected subgroup whose semi-simple part coincides with
the
 semisimple part of the centralizer of
a toral subgroup of $G$, such that the coset space $G/H$ is
even-dimensional. Then $G/H$ has a homogeneous complex structure
and each C-space is homeomorphic to such coset.
\end{theorem}

Let ${\mathfrak g}^c$ be a complex semisimple Lie algebra and $G^c$ a complex Lie group with Lie algebra ${\mathfrak g}^c$. Wang has shown that $G^c$ can be chosen so that $G\subset G^c$ and there is complex Lie subgroup $H^c\subset G^c$ such that $G/H = G^c/H^c$ and $H=H^c\cap G$.

 In order to detect
the spaces with vanishing first Chern class we need some basic
facts about semisimple Lie algebras and their parabolic
subalgebras.
  Fix a Cartan subalgebra
  (maximal toral subalgebra) $\mathfrak t$ in ${\mathfrak g}^c$. Then we have
a system of roots $R \subset {\mathfrak t}^*$ defined by $\mathfrak t$ in ${\mathfrak g}^c$. There is also a
 distinguished set of simple roots $\Pi$
 in $R$ which forms a basis for ${\mathfrak t}^*$ as a (complex) vector
 space and defines a splitting $R=R^+\bigcup R^-$ of $R$ into positive and negative roots. Let
 ${\mathfrak g}$ and ${\mathfrak h}$ are the (real) Lie algebras of $G$ and $H$. The above mentioned article of Wang provides the inclusion ${\mathfrak h}^c_{ss} \subset
{\mathfrak h}^c \subset
 {\mathfrak j}^c$, where ${\mathfrak j}^c$ is a parabolic subalgebra, which
is a centralizer of a torus and ${\mathfrak j}^c_{ss} = {\mathfrak h}^c_{ss}$. Here the subscript
 $"ss"$ denotes the semisimple part. In particular

$$ {\mathfrak h}^c = {\mathfrak a} + {\mathfrak h}^c_{ss}$$
where $\mathfrak a$ is a commutative subalgebra of the Cartan subalgebra $\mathfrak t$ of ${\mathfrak g}^c$. The
parabolic algebra ${\mathfrak j}^c$ satisfies ${\mathfrak j}^c = {\mathfrak t} + {\mathfrak j}^c_{ss}$ and is equal to
the normalizer of ${\mathfrak h}^c$ in ${\mathfrak g}^c$. The sum here is not direct because part of ${\mathfrak
t}$ is contained in ${\mathfrak j}^c_{ss}$. Let $J^c$ be a parabolic subgroup of $G^c$ with algebra ${\mathfrak
j}^c$ and $G^c/J^c$ is the corresponding (generalized) flag manifold (called also rational homogeneous space).
Then according to Wang \cite{W}, there is a holomorphic fibration $G^c/H^c\rightarrow G^c/J^c$ with fiber a
complex torus, which determines the complex structure on $G/H$.

 As is shown by Tits \cite{T}, for general complex Lie group $G^c$ this
 fibration is canonically defined and unique in a sense that any other
homogeneous fibration over $G^c/H^c$ factors through $G^c/J^c$. It is
called the {\it Tits
 fibration}. Moreover he shows that the fiber is a complex parallelizable manifold and
in the case of compact and semisimple group $G$ the fibration is a principal toric bundle. This reduces the study of $G/H$ to a study of the toric
fibrations over the flag manifold $G/J$, where $J=J^c\cap G$. Moreover any compact Lie group can be represented as $G=G_1\times T^k/H$ for a semisimple $G_1$ and a torus $T^K$ where $H$ is finite abelian subgroup of the center of $G$. Then, up to a finite abelian cover, the Tits fibration of $G/H$ for a compact $G$ is also a principal toric bundle.

 Now every complex structure on the flag manifold is determined by an ordering of the
system of roots of ${\mathfrak g}^c$.  The complex structure on $G/J$ defines subset $\Pi_0$ in $\Pi$ which
corresponds to ${\mathfrak j}^c$. Since this correspondence determines the second cohomology and the first Chern
class  of $G/H$ we provide more details about it.
  In general ${\mathfrak j}^c_{ss}$ is determined
by the span of all roots in $R$ which are positive with respect to
$\Pi_0$. Then the complement $\Pi - \Pi_0 = \Pi'$ provides a basis
for
 the center $\zeta$ of ${\mathfrak j}^c$ and there is an
identification $span_{\bf Z}(\Pi') = H^2(G/J,\mathbb{Z})$. The
identification \cite{Al} is:

$$ \xi \rightarrow \frac{i}{2\pi} d\xi$$
where $\xi$ is considered as a left invariant 1-form on G which is a subgroup of $G^c$ and $d\xi$ is
$ad({\mathfrak j})$-invariant, hence defines a 2-form on $G/J$. This form is obviously closed and in fact
defines non-zero element in $H^2(G/J,\mathbb{Z})$. Moreover every class in $H^2(G/J,\mathbb{Z})$ has unique
representative of this form.

 Now we are interested in the first Chern class of G/H. It is determined by the
so-called Koszul form for $G/H$ (see \cite{T},\cite{AlP}). The definition of
this form is

$$ \sigma_{G/H}(X) = Tr_{\frac{\mathfrak g}{\mathfrak h} } (ad(IX)-Iad(X)), X \in {\mathfrak g}$$
where $I$ is the invariant complex structure on $G/H$, extended as $0$ on ${\mathfrak h}$. Note that the form itself is defined only on  $G$ and without the restriction of $G$ being semisimple. According to
\cite{AlP}, Proposition 4.1 for semisimple Lie groups $G$,
 \[
 \sigma_{G/H}=2i(\sigma_G -
\sigma_H)
\]
 where $\sigma_G$ is the sum of positive roots in ${\mathfrak g}^c$
and $\sigma_H$ is the sum of positive roots in ${\mathfrak g}^c$ which are also in ${\mathfrak h}^c$. From here
 one has that $\sigma_{G/J}
 = \sigma_{G/H}$ for the Tits fibration $\pi: G/H
 \rightarrow G/J$  since the semisimple parts of ${\mathfrak j}^c$ and ${\mathfrak h}^c$ coincide.
  As is proved by Koszul, the form $d\sigma$ descends to $G/H$ and
 represents the first
Chern class of the homogeneous manifold $G/H$. The same is true for $G/J$.

\begin{theorem}\label{c1=0} Let $G$ be a (connected) compact semisimple Lie group and $H$ is a closed subgroup, such that $G/H$ has a $G-$invariant complex structure $I$. The first Chern class of $(G/H, I)$ vanishes iff
$$\sigma_{G/H}|_{\mathfrak a} = 0$$
i.e. the restriction of $\sigma$ to $\mathfrak a$ vanishes.
\end{theorem}

\bproof The key point is that $d\sigma$
 defines the zero element in $H^2(G/H,\mathbb{Z})$  iff $\sigma_{G/H}$
descends to an invariant 1-form on $G/H$ itself. But if $d\sigma_{G/H}=df$ for some 1-form $f$ on $G/H$, than one
can symmetrize $f$ as in the section 3.2 below to obtain an invariant form $f'$ on $G/H$ with $df=df'$. Then the
pull-back of $f'$ is an invariant form on $G$ with the same differential as $\sigma_{G/H}$, hence it coincides
with $\sigma_{G/H}$.  As it is proved in Tits (see \cite{T},\cite{Al}) $\sigma_{G/H}$ is a sum with positive
integer coefficients of elements of $\Pi'$ and thus is in $\zeta^*$, the dual of the center $\zeta$ of ${\mathfrak
j}^c$. Then it descends to $G/H$ iff it vanishes on ${\mathfrak a}$. \eproof

With this characterization in mind we continue with some examples.


\subsection{Examples of C-spaces with vanishing first Chern class} We concentrate on quotients of a compact simple classical Lie group,
although the Calabi-Eckmann manifolds show that some homogeneous
manifolds may be reducible as smooth manifolds and irreducible when considered
as complex homogeneous ones.

 The first type of  examples include the following two extreme cases:

$i$) ${\mathfrak a}= 0$ i.e. $\mathfrak h$ is semisimple itself.

$ii$) $H = U(1)$, where U(1) is appropriately embedded in odd-dimensional $G$, so that the principal bundle $G
\rightarrow G/U(1)$ is a pull back of the canonical $U(1)$-bundle over the full flag manifold of $G$.

For the first case we start with an example from the $A_{\ell}$-series. Consider
\[
M = SU(n)/SU(n_1)\times SU(n_2)\times ... \times SU(n_k), k - odd, n_i >1
\]
Here $SU(n_1)\times ...\times SU(n_k)$ is diagonally embedded as a
matrix group in $SU(n)$
 and $n_1 + n_2 + ... + n_k = n$.
 The condition $k - odd$ ensures that the space is even dimensional. Then the Tits
 fibration is $M \rightarrow
SU(n)/S(U(n_1)\times U(n_2)...\times U(n_k))$ with fiber $T^{k-1}$. The existence of a complex structure follows
by Theorem \ref{c1=0}. The vanishing of the first Chern class follows by the above theorem. These examples clearly
could be generalized to
\[
 M =SU(n)/SU(n_1)\times SU(n_2)\times ... \times SU(n_k), n_1 + n_2 ... + n_k \leq n,
n_i >1 \]
 where $n - (n_1 + n_2 + ... + n_k) + k$ is odd.

 For the other classical compact Lie groups we have also the following spaces:
\[
\begin{array}{llll}
M &= &SO(2n)/SU(n_1)\times ... \times SU(n_{2k})\times SO(2l),& n_1 +...+n_{2k}+l=n \\
M& = & SO(2n)/SU(n_1)\times ... \times SU(n_{2k})\times SU(n_{2k+1}), & n_1 +...+n_{2k+1}=n \\
M&=& SO(2n+1)/SU(n_1)\times ...\times SU(n_{2k})\times SO(2l+1),&
n_1 +...+n_{2k}+l
=n \\
M&=& Sp(n)/SU(n_1)\times ... SU(n_{2k})\times Sp(l), & n_1 +...+n_{2k}+l=n \\
\end{array}
\]
We assume above that not all of $n_i$ vanish.

{\it Case ii:} The example is $G/U(1)$, but with U(1) appropriately embedded. We concentrate on $G=SU(n)$. The Tits
fibration in this case is $SU(n)\rightarrow F_{1,1...1}$ where $F_{1,1...,1}$ is the standard flag manifold
$SU(n)/S(U(1)\times ...\times U(1))$. In this case we can avoid more complicated tools like the painted Dynkin
diagram bellow because $\sigma_{F} = \sigma_{SU(n)} - \sigma_{S(U(1)\times...\times U(1))} = \sigma_{SU(n)}$,
where $\sigma_G$ is the sum of all positive roots of $G^c$. Then $\sigma_{S(U(1)\times...\times U(1))} = 0$ because $S(U(1)\times...\times U(1))$ is abelian. To describe the sum of positive roots we need some notations which appear in \cite{FH}. Let
$L_i$ be the matrix with $i$-th diagonal element equal to 1 and all others being 0. Then the set of all roots of
the complexified Lie algebra $sl(n,\mathbb{C})$ is $e_{i,j} = e_i - e_j$, where $e_i$ are the duals of $L_i$. A
set of simple roots is $e_{i,i+1}$ which also determines the positive roots $e_{i,j}, i<j$. Then the sum of all
positive roots is:
\[
 \sum_{i<j}e_{i,j} = \sum_{k=1}^{n-1} k(n-k)e_{k,k+1} = \sum_{k=1}^{n}(n-2k+1)e_k
\]
Then we have:
\begin{proposition}\label{example1}
The space $SU(n)/U(1)$ for $n$ even endowed with the homogeneous complex structure from \cite{W} has
vanishing first Chern class iff $U(1)$ is embedded as a set of
diagonal matrices:

$$A=diag(e^{2\pi\theta_1 t},e^{2\pi\theta_2 t},...,e^{2\pi\theta_n
t})$$
 with $\theta_n=-\theta_1-...-\theta_{n-1}$ and satsfying
 $$\sum_{k=1}^{n}(n-2k+1)\theta_k = 2\sum_{k=1}^{n-1}(n-k)\theta_k=0$$
\end{proposition}

Next we consider the general case of factors of $SU(n)$. In this
case the Tits fibration is of the form:
$$ SU(n)/SU(n_1)\times...\times SU(n_k)\times T^l \rightarrow
SU(n)/SU(n_1)\times...\times SU(n_k)\times T^m$$ with $n_1 + ...+n_k + m - k= n-1$ and $l<m$. At this point we
notice that an invariant complex structure on the flag manifold $SU(n)/SU(n_1)\times...\times SU(n_k)\times T^m$
is not unique and depends on the so called {\it painted Dynkin diagram} (or black-white Dynkin diagram). The
painted diagrams are used to describe (generalized) flag manifolds and homogeneous Einstein metrics on them
(see for example \cite{Andr}). For a flag manifold, painted Dynkin diagram is obtained by blackenning the vertices which
correspond to $\Pi'$. We refer to \cite{Al} for the details and use directly their result for two particular
examples which enlighten the general case.

{\it Example of $A_l$-type.} Choose the flag manifold for the base of the Tits fibration to be
$SU(11)/SU(4)\times SU(3)\times SU(2)\times T^4 = SU(11)/S(T^2\times U(4)\times U(3)\times U(2))$. It
corresponds to
 a painted Dynkin diagram:

\begin{tabular}{cl}
\hspace{1cm} &\AnblackMarsOK\\
\end{tabular}
\vspace{.2in}

The diagram also determines the complex structure on the flag manifold above. Now we need the Koszul
form $\sigma$ of this flag manifold. In case of the $A_n$-series it is described in \cite{Al}, Propositrion 4.1
and Proposition 5.2. From there we have that
$$\sigma = (2+b_1)\overline{\alpha_1} +
(2+b_2)\overline{\alpha_2}+ ... + (2+b_m)\overline{\alpha_m}$$
where $\overline{\alpha_i}$ are the fundamental weights
corresponding to the roots with black circles. They are defined as
\begin{equation}\label{fundweights}
\frac{(\overline{\alpha_k},\alpha_j)}{(\alpha_j,\alpha_j)} =
\delta^i_j
\end{equation}
where $(,)$ is the product arising from the Killing form.
 The numbers $b_i$ are nonnegative and in our particular case are equal to the number of
white circles of the Dynkin diagram, which are connected with the
black circle corresponding to the root $\alpha_i$ by a series of
white circles \cite{Al}.

In particular for the above diagram we have:
$$ \sigma = (2+b_1)\overline{e_{1,2}} +
(2+b_2)\overline{e_{2,3}} + (2+b_3)\overline{e_{6,7}} +
(2+b_4)\overline{e_{9,10}} = 2\overline{e_{1,2}} +
5\overline{e_{2,3}}+ 7\overline{e_{6.7}}+5\overline{e_{9,10}} $$

Now to obtain explicit expression of the above element $\sigma$ in terms of $e_{i,i+1}$ we need a description of
the fundamental weights. This could be done following \cite{FH}.

 On $\mathfrak{ g}^c=sl(11,\mathbb{C})$ the matrices $e_i$ are
orthonormal with respect to $(,)$. Then one can check directly that the following elements satisfy the condition
(\ref{fundweights})
$$
L_k = e_1 +...+e_k - k/n(e_1+...+e_n) = \overline{e_{k,k+1}}
$$
So in our case we have

\begin{eqnarray}
\nonumber
\sigma &=&2e_1-2/11(e_1+...+e_{11})+5(e_1+e_2)-10/11(e_1+...+e_{11})+7(e_1+...+e_6)\\
\nonumber
 & &-42/11(e_1+...+e_{11})+5(e_1+...e_9)-45/11(e_1+...+e_{11})\\
\nonumber
 &=&10e_1+8e_2+3(e_3+e_4+e_5+e_6)-4(e_7+...+e_9)-9(e_{10}+e_{11})\\
\nonumber
\end{eqnarray}

The dimension count gives $dim SU(11)/S(T^2\times U(4)\times U(3)\times U(2)) = 90$ and there could be either a 2 or
4-dimensional torus as a fiber for the Tits fibration with this base. We also consider the subgroup $J=S(T^2\times
U(4)\times U(3)\times U(2))$ embedded in $SU(11)$ in the standard diagonal form with unitary blocks of order
1,1,4,3, and 2 respectively so that $G/J$ is generalized flag manifold. The case of 4-dimensional fiber leads to an example of the previous type because the
subgroup $H$ will be semisimple. So we consider the 2-dimensional fibers. At the Lie algebra level we have to
add appropriate 2-dimensional space $\mathfrak a$ of diagonal matrices to the Lie algebra ${\mathfrak j}_{ss}$, the real part of $\mathfrak{j}^c_{ss}$.
It has to be of the form
$$ {\mathfrak a} =
diag(x_1,x_2,x_3,x_3,x_3,x_3,x_4,x_4,x_4,x_5,x_5)
$$
and should obey the following conditions:
\begin{eqnarray}
\nonumber
x_1+x_2+4x_3+4x_4+2x_5 &=& 0 \\
\nonumber
10x_1+8x_2+12x_3-12x_4-18x_5 &=& 0 \\
\nonumber
\end{eqnarray}
The first equation comes form the requirement that the matrices in $\mathfrak a$ are trace-free. The second
follows from Theorem \ref{c1=0} and the form of $\sigma$ above. Now we can fix two linearly independent integer
solutions $(v_1,...,v_5)$ and $(w_1,...,w_5)$ of these equations. Then the Lie algebra ${\mathfrak h}=
{\mathfrak j}_{ss}+{\mathfrak a}$ should be:
$$
{\mathfrak h} = diag(v_1t+w_s,v_2t+w_2s,
(v_3t+w_3s)A,(v_4t+w_4s)B, (v_5t+w_5s)C)
$$
where $A,B,C$ are trace-free skew-adjoint matrices of order 4,3
and 2 respectively. Then at the end we obtain that $SU(11)/H$ is a
complex homogeneous manifold with vanishing first Chern class, if
$H$ is of the form:
$$
H =
diag(e^{2i\pi(v_1t+w_1s)},e^{2i\pi(v_2t+w_2s)},e^{2i\pi(v_3t+w_3s)}A,
e^{2i\pi(v_4t+w_4s)}B, e^{2i\pi(v_2t+w_2s)}C)
$$
where $A,B,C$ are unitary matrices with determinants one and order
4,3 and 2 respectively. Moreover any such manifold with a
stationary subgroup $H$ containing strictly $Id_2\times
SU(4)\times SU(3) \times SU(2))$ is of this form.

{\it Example of $C_l$-type.} As a last example we consider the following quotient of the symplectic group  -
$\frac{Sp(7)}{SU(2)\times Sp(2)\times T^2}$. Let $sp(2n, \mathbb{C})$ be the complex Lie algebra of symplectic
complex $2n\times 2n$ matrices. Its diagonal has a basis of the form $L_i = E_{i,i}-E_{n+i,n+i}$, where
$E_{i,j}$ is the matrix with $(i,j)$-th element 1 and all others 0. Then a set of simple roots is given by the
duals to the matrices $L_i-L_{i+1}$ and $2L_n$. Denote by $e_i$ the dual to $L_i$ and then the simple roots are given
by $\alpha_i=e_i-e_{i+1}, \alpha_n=2e_n$. To calculate the Koszul form we can not use the routine from the
previous example but we can adapt the algorithm from \cite{BFR}. Bellow is the Dynkin diagram for $sp(14,
\mathbb{C})$ together with the coefficients for each simple root, with which it enters in the sum of all
positive roots:

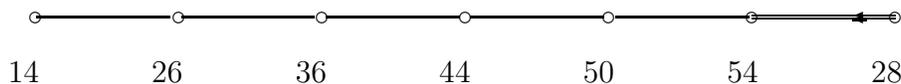
\begin{tabular}{cl}

\hspace{0.65cm}  & \CnMarsOK\\

\end{tabular}

$$
14\hspace{1.5cm}26\hspace{1.5cm}36\hspace{1.5cm}44\hspace{1.5cm}50\hspace{1.5cm}54\hspace{1.5cm}28
$$

 The similar diagram corresponding to $sl(2,\mathbb{C})\oplus sp(4,\mathbb{C})$ is

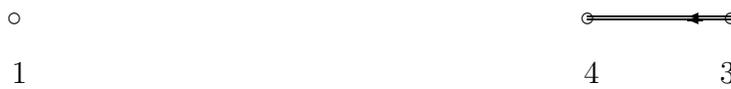
\begin{tabular}{cl}

\hspace{0.3cm}  & \CnwhiteMarsOK\\

\end{tabular}

 $$
 \hspace{1.5cm}1\hspace{7.4cm} 4\hspace{1.6cm}3
 $$

Now we take the difference of the two diagrams: Complete the second diagram with black dots corresponding to the
complementary roots and take the difference of the corresponding coefficients. The result is given by (notice
that our painting is opposite to the one in \cite{BFR}):

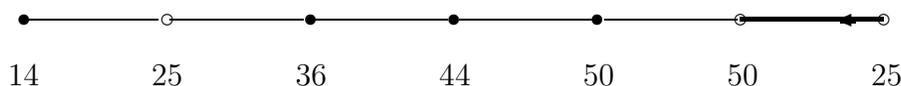
\begin{tabular}{cl}

\hspace{.5cm}  & \CnblackMarsOK\\

\end{tabular}

$$
14 \hspace{1.5cm}25\hspace{1.5cm}36\hspace{1.5cm}44\hspace{1.5cm}50\hspace{1.5cm}50\hspace{1.5cm}25
$$
This means that the Koszul form of the flag manifold $\frac{Sp(7)}{SU(2)\times Sp(2)\times T^4}$ is given by
$$\sigma = 14\alpha_1+25\alpha_2+36\alpha_3+44\alpha_4+50(\alpha_5+\alpha_6)+24\alpha_7 =
14e_1+11(e_2+e_3)+8e_4+6e_5$$

Now we consider the compact real forms of the Lie algebras above. The compact form of $sp(2n, \mathbb{C})$ is denoted
simply by $sp(n)$ and is the the set $n\times n$ quaternionic matrices $A$ with $A+\overline{A} =0$. Then the
$sp(2)$ which corresponds to the last two dots in the painted diagram above is represented by a $2\times 2$
submatrix in the lower right corner. In the same way $su(2)$ is represented by one diagonal entry which is an
imaginery quaternion. The 2-dimensional abelian subalgebra which should be added to the sum $su(2)+sp(2)$ is the
algebra of diagonal matrices with entries $i(\lambda_1t+\mu_1s,\lambda_2t+\mu_2s, \lambda_2t+\mu_2s,
\lambda_3t+\mu_3s,\lambda_4t+\mu_4s,0,0)$, for which ${\bf \lambda}=(\lambda_1,\lambda_2,\lambda_3,\lambda_4)$
and ${\bf \mu}=(\mu_1, \mu_2, \mu_3, \mu_4)$ are independent integral solutions of $14x+22y+8z+6v=0$. Then by
Theorem \ref{c1=0} we have that $M=\frac{Sp(7)}{H}$ is a complex homogeneous space of vanishing first Chern
class if $SU(2)\times Sp(2)\times T^2 = H$ is embedded as $H=diag(e^{2\pi i(\lambda_1t+\mu_1s)}, q, e^{2\pi
i(\lambda_2t+\mu_2s)}, e^{2\pi i(\lambda_3t+\mu_3s)}, e^{2\pi i(\lambda_4t+\mu_4s)}, A)$, where $A$ runs over
$Sp(2)$ and $q$ over $SU(2)=Sp(1)$.
\begin{remark} If the integer solutions in the above two examples are chosen with greatest common divisor one, than the factor-spaces are simply-connected. The topology of the examples from Proposition \ref{example1} is discussed in Section 7.
\end{remark}
\section{CYT metrics on C-spaces}
\subsection{Canonical connections and CYT metrics}
\noindent We provide here the properties of the canonical connections on Hermitian manifolds which we will need. Let $(M,J,g)$ be Hermitian manifold and $F(X,Y)=g(JX,Y)$ is the K\"ahler form.  Let $d^c$ be the operator $(-1)^nJdJ$ on n-forms and $D$ be the
Levi-Civita connection of the metric $g$. Then, a family of canonical connections is given by \cite{G1}:
\begin{equation}
g(\nabla^t_AB, C)
=g(D_AB, C)+\frac{t-1}{4}(d^cF)(A,B,C)+\frac{t+1}{4}(d^cF)(A, JB, JC),
\end{equation}
where $A, B, C$ are any smooth vector fields and the real number $t$ is a free parameter.
The connections $\nabla^t$ are called \it canonical \rm connections.
It is known that the connection $\nabla^{1}$ is the Chern connection on the holomorphic
tangent bundle. The connection $\nabla^{-1}$ is called KT connection by physicists and Bismut connection by some mathematicians. Its mathematical features and background are articulated in \cite{G1}. It is apparent in the above
expression that when the Hermitian structure is a K\"ahler structure, this one-parameter family
of canonical connections collapses to a single connection: the Levi-Civita connection.

 From formula (2.7.6) in \cite{G1} we have
$$
\nabla^t s - \nabla^u s = i \frac{t-u}{2} \delta F \otimes s
$$
for any section $s$ of the canonical bundle $K^{-1}$ and $\delta$ is the codifferential.

Let $R^{t}$ be the curvature of $\nabla^t$ and $\rho^t(X,Y)=\Sigma g(R^t  (X,Y)E_i,JE_i)$ be the corresponding
trace. Then $i\rho^t$ is the curvature of $K^{-1}$ and from the above relation we obtain:
$$
\rho^t - \rho^u = \frac{t-u}{2} d\delta F
$$
Denote by $\rho$ and $\rho^B$ the
Ricci forms of the Chern and KT connections respectively.







The formula (53) from \cite{G2}  gives representation  in local coordinates:
\begin{equation}\label{chern}
\rho = -1/2 dd^c\log(det(g_{\alpha\overline{\beta} }))
\end{equation}
in particular $\rho$ is always a (1,1)-form. In general $d\delta F$ is not necessary a (1,1)-form.

\begin{definition}
  A Hermitian metric is called {\it CYT} (or Calabi-Yau with torsion)
if $\rho^B = 0$. Then $(g,J)$ is called a CYT structure and $\nabla^B$ a CYT connection. Equivalently, the (restricted) holonomy group of the Hermitian connection with totally skew-symmetric torsion is contained in $SU(n)$.
\end{definition}

  For later use we recall also the toric bundle construction of CYT structures described in \cite{GGP}. Suppose that
$X\rightarrow B$ is a $T^{2k}$-bundle over complex manifold $B$ with characteristic classes
$\omega_1,...,\omega_{2k}$. If $\eta_1,...,\eta_{2k}$ are connection 1-forms with $d\eta_i=\omega_i$, then there
is an almost complex structure on $X$ defined in the following way. On horizontal cotangent subspaces it is the
horizontal lift of the structure on the base $B$. The vertical subspaces are spanned by $Ker(\eta_i)$ and the
definition is just $J\eta_{2i-1}=\eta_{2i}$. The result is an almost complex structure $J$ on $X$ which is
integrable if $\omega_{2i-1}+i\omega_{2i}$ is a 2-form of type (1,1)+(2,0). Moreover any toric invariant complex
structure on $X$ for which the bundle projection is holomorphic arises  in this way for some choice of the
connection.

 Suppose now that
$g_X$ is a Hermitian metric on the base manifold $X$
with K\"ahler form $F_X$. Consider a Hermitian metric on $M$
defined by
\[
g_M = \pi^*g_X + \sum_{\ell=1}^{2k} \theta_\ell^2.
\]
Since $J\theta_{2j-1}=\theta_{2j}$, the  K\"ahler form for this Hermitian metric
is
\[
F_M= \pi^*F_X+ \sum_{j=1}^k\theta_{2j-1}\wedge \theta_{2j}.
\]

By \cite{GGP} the Ricci forms of the canonical connections of $X$ and $M$ are related by:
\[
\rho^t_M=\pi^*\rho^t_X+\frac{t-1}2\sum_{\ell=1}^{2k}d((\Lambda\omega_\ell) \theta_\ell).
\]
where $\Lambda\omega_\ell=g(F_X,\omega_\ell)$. In particular when $g(F,\omega_\ell)=constant$ we have
\begin{equation}\label{torusrel}
\rho^B_M=\pi^*(\rho^B_X-\sum_{\ell=1}^{2k}g(F_X,\omega_\ell)\omega_\ell).
\end{equation}
For more details see \cite{GGP, GP}.

\subsection{Existence of CYT metrics on C-spaces}

Many examples of CYT spaces are provided by the compact homogeneous C-spaces of Wang. Assume that $M= G^c/L$ is a
compact complex homogeneous space defined by a transitive action of a complex Lie group $G^c$ with finite
fundamental group. Then by an observation in \cite{W}, Remark (2.3), $M$ is analytically isomorphic to
$G/(G\cap L)$, where $G$ is a maximal connected compact semisimple subgroup of $G^c$. For general compact Lie group actions we have:

\begin{theorem}\label{compact} Let $G$ be a (connected) compact Lie group and $H$ a closed subgroup, such that the homogeneous space $M=G/H$ admits a $G$-invariant complex structure with vanishing first Chern class. Then there is a 1-parameter family of invariant complex structures which contains the given structure and a structure that admits a compatible CYT metric.
\end{theorem}

\bproof
  Consider the splitting of the Lie algebra ${\mathfrak g} = {\mathfrak h} + {\mathfrak m}$
of $G$, where ${\mathfrak m}$ is the orthogonal compliment of $\mathfrak h$ in ${\mathfrak g}$ with respect to the metric $B(X,Y)$
given by the negative of the Killing form on $G$. Now $B$ is bi-invariant and in particular Ad$G$-invariant
and when restricted to ${\bf h}$ and ${\bf m}$ is non-degenerate. Then by theorem 3.4, Chapter 10 of \cite{KN}, we have that
$M$ is a naturally reductive homogeneous space. Such spaces have a canonical homogeneous connection which is
metric and its torsion at the identity coset is given by $T(X,Y)_0 = -[X,Y]|_{\mathfrak m}$. Then for naturally
reductive spaces $B(X,[Y,Z]|_{\bf m})$ is totally skew-symmetric 3-form. Suppose for the moment that the complex structure is Hermitian with respect to the metric induced by $B$. Now since any invariant tensor field is parallel with respect to the
canonical connection, it is the Bismut connection for the naturally reductive complex homogeneous space $M$. We will show
 that it is CYT. We recall first that the canonical bundle $K$ for the complex structure is topologically trivial
because its first Chern class vanishes. Then we can find a (n,0)-form which is nowhere zero and by
symmetrization, we obtain a global $G$-invariant (n,0)-form. More precisely if $\omega$ is non-vanishing
(n,0)-form we define:
$$
\mu(\omega)_x(A_1,...,A_n) = \int_G \omega_{g(x)}(L_g^*(A_1),...,L_g^*(A_n)) dm_g
$$
where $dm_g$ is bi-invariant volume form on $G$. Then $\mu(\omega)$ is invariant. Since the invariant tensors
are parallel, it is parallel with respect to the canonical connection. Then the only point is to prove that it
is everywhere non-zero. However $\mu(\omega)\wedge\overline{\mu(\omega)} = f dm_g$ for some function $f$. Since
$\mu(\omega)$, $\overline{\mu(\omega)}$ and $dm_g$ are $G$-invariant, $f=constant$. So $\mu(\omega)$ is either
nowhere zero or vanishing identically. However in the definition of $\omega$ we have the freedom to multiply it
by arbitrary nonvanishing real function. Assume that $\mu(\omega)=0$ for any choice of $\omega$. Fix tangent
(1,0)-vectors $(A_1,...,A_n)$ at the identity coset in $G/H$. Consider the function $h(g)=
L_g^*\omega(A_1,...,A_n)$. Then our assumption implies that $\int_G f(g)h(g)dm_g=0$ for every smooth function $f(g)$
such that $f(gh_1)=f(gh_2)$ for every $h_1,h_2\in H$ and $f(g)\neq 0$. By taking limit, we can actually extend
our assumption to every integrable function $f$. Since $h(e)\neq 0$, we can assume that its real part is $\Re
h(e)>0$ in some neighborhood of $H$ in $G$. Then by taking cut-off real function $f$ we can restrict the above
integral on a smaller neighborhood of $H$ and see that its real part is strictly positive, which is a contradiction.
Then the proof follows from the next Lemma. \eproof

\begin{lemma}
If $G$ is compact then every $G$-invariant complex structure on $G/H$ is either compatible with $B$ or belongs to a 1-parameter family of homogeneous complex structures which contains a structure compatible with $B$
\end{lemma}

\bproof
Since any invariant complex structure admits a Tits fibration, it is determined by an invariant complex structure on the base generalized flag manifold and a complex structure on the torus fiber (possibly after taking a finite cover). The main observation is that any invariant complex structure on the base is compatible with the induced $B$. A choice of $B$-orthonormal invariant 1-forms which provide a basis for $\mathfrak{a}^*$ will determine a connection 1-forms of a "canonical" connection for the Tits fibration. Then the given complex structure on the fiber is determined by a constant matrix $J$ of square minus identity in such basis. The structure is $B$-compatible if and only if this matrix is orthogonal. However if it is not, we can find a path of such matrices which connects $J$ to an orthogonal one, because the set $O(2n)/U(n)$ is a deformation retract of $GL(2n,\mathbb{R})/GL(n,\mathbb(C)$. Such a path will define a 1-parameter family an invariant almost complex structures on $G/H$. Since the second cohomology group of a generalized flag manifold contains only classes of type $(1,1)$, characteristic classes of the Tits fibration are of type $(1,1)$. Then all almost complex structures in the family are integrable by the results in \cite{GGP}.
\eproof

\begin{remark} In Proposition 3 below we show that there are $G$-invariant complex structures which do not admit any $G$-invariant CYT metric on the group $G=SU(2)\times SU(2)$. This shows that the deformation above is necessary. It is not clear whether there is any CYT metric in the general case.
\end{remark}

\subsection{Deformations of the bi-invariant metric on SU(3) and CYT structures on $S^3\times S^3$}

Suppose that we have a complex semisimple Lie algebra ${\mathfrak g}^c$ and a root decomposition with root basis
$(E_{\alpha},E_{-\alpha}, H_k)$, where $H_k$ is a basis of a fixed Cartan subalgebra $\mathfrak{k}$. Suppose its compact
form has associated Lie group $G$ which admits complex structure $I$ with $I(E_{\alpha})= iE_{\alpha}, I(E_{-\alpha})=-iE_{-\alpha}$. Then
consider a left invariant Hermitian metric on $G$ such that it is also right-invariant for the maximal tori. In
such case it is of the form $ g(X,Y)= B(\Lambda(X),Y)$ where $B(.,.)$ is the Hermitian extension of the bi-invariant metric on $G$ and $\Lambda$ is
a positive Hermitian operator defined by $\Lambda(E_{\alpha})=\lambda_{\alpha} E_{\alpha}$ for $\alpha>0$ and any
appropriate matrix on $H_k$. Consider the case when it is identity for $H_k$. Then we calculate $\delta F$ from the identity:
$$\delta F\circ I = 1/2F\lrcorner dF$$

Because $dF(X,Y,Z)=F([X,Y],Z)+F([Y,Z],X)+F([Z,X],Y)$ for invariant metrics, the sum $\sum_i dF(E_i,IE_i,Z)$ has
contributions only from $dF(E_{\alpha},E_{-\alpha}, H) =$ $ -\lambda_{\alpha} \alpha(IH)$, for $H\in \mathfrak{k}$. Then $$\delta F\circ
I = - 1/2 \sum_{\alpha} \lambda_{\alpha}.\alpha \circ I$$ where the sum is over all positive roots. By Koszul \cite{K}
the Ricci form of the Chern connection for $I$ is the exterior derivative of the half-sum of all positive roots
for any invariant metric. So the Ricci form $\rho^B$ is given by
$$\rho^B = 1/2\sum_{\alpha}(1-\lambda_{\alpha})d\alpha$$
for the sum running over all positive roots again.

In the case of $SU(3)$, the roots are $(\alpha,\beta,\alpha+\beta
,-\alpha,-\beta,-\alpha-\beta)$.Then the sum is
$2\alpha+2\beta$ and we have:
\[
\rho^B=0 \hspace{.1in}iff \hspace{.1in} \lambda_{\alpha} + \lambda_{\alpha+\beta} =
\lambda_{\beta}+\lambda_{\alpha + \beta} = 2
\]

Then taking $\lambda_{\alpha}$ as a parameter, we have that:

\begin{proposition} For a fixed invariant complex structure on $SU(3)$, there is a 1-parameter family of invariant CYT metrics parametrized by the interval (0,2).
\end{proposition}

One can compare the deformation of the metric above with the fact that any small hypercomplex deformation of the HKT structure on $SU(3)$ is HKT \cite{GPP}. Since every HKT structure is CYT, this implies that for many non-homogeneous complex structures on $SU(3)$ there are CYT metrics.

\vspace{.1in}
{\it CYT structures on $SU(2)\times SU(2)$}.
The space $SU(2)\times SU(2)$ carries a Calabi-Eckmann complex structure and has a basis of global left-invariant 1-forms $\alpha_i,e_i^+,e_i^-$ where $i=1,2$ correspond to the two factors. It satisfies the Maurer-Cartan equations $d\alpha_i=e_i^-\wedge e_i^+, de_i^+=\alpha_i\wedge e_i^+, de_i^-=-\alpha_i\wedge e_i^-$. The forms $\alpha_1,\alpha_2$ define the canonical connection for the Tits fibration which in this case is $SU(2)\times SU(2)\rightarrow {\bf CP}^1\times {\bf CP}^1$. The forms $e_i^{\pm}$ are pull-backs from forms on the base. Every left-invariant complex structure on $SU(2)\times SU(2)$ is also right-invariant with respect to the torus fibers of the Tits fibration. Then, up to sign on the factors of the base $ {\bf CP}^1\times {\bf CP}^1$, it is determined by the following:
\begin{equation}\label{cs}
J(\alpha_1)=a\alpha_1+b\alpha_2, J(\alpha_2)=c\alpha_1+d\alpha_2, J(e_i^+)=e_i^-, J(e_i^-)=-e_i^+
\end{equation}
where $\left(\begin{array}{ll}
 a&b\\
c&d
\end{array}\right)^2=-Id$, so
$d=-a, b\neq 0, c=-\frac{a^2+1}{b}$. From the previous sections we know that the Ricci form of the Chern connection for any invariant Hermitian metric is given by $d\sigma_{SU(2)\times SU(2)} = d(\alpha_1+\alpha_2)$. If the K\"ahler form of the metric is $F$ then from $0=\rho^B=\rho^{Ch}-d\delta F$ we have $\rho^{Ch}=d\delta F$ for the equation of a CYT metric. Since $\rho^{Ch}=d\sigma_{SU(2)\times SU(2)}$ and both $\sigma_{SU(2)\times SU(2)}$ and $\delta F$ are invariant 1-forms then  $\sigma_{SU(2)\times SU(2)}=\delta F + \beta$ for an invariant and closed 1-form $\beta$. However $\beta$ is not exact and since $SU(2)\times SU(2)$ is simply-connected, $\beta=0$. This considerations are obviously valid for any compact semisimple Lie group $G$, so we obtain:
 \begin{lemma} Let $G$ be a compact semisimple Lie group with a fixed left-invariant complex structure on it. Then
 the CYT equation for any compatible left-invariant metric $g$
  can be written as:
  \begin{equation}\label{compactG}
  \sigma_G=\delta F
  \end{equation}
  where $F$ is the K\"ahler form of $g$ and $\delta$ is the codifferential.
  \end{lemma}

To solve this equation on $SU(2)\times SU(2)$ we denote by $g$ the metric induced on the cotangent bundle and notice that
(\ref{compactG}) is equivalent to $g(\alpha_1+\alpha_2, \xi_i)=g(\delta F, \xi_i)=g(F,d\xi_i)$ for any basis $\xi_i$ of left-invariant 1-forms. The last equation follows by integration by parts and the fact that all functions involved in the equalities are constants. We choose the basis $\xi_i$ to be $\alpha_i, e_i^{\pm}$. Using the Maurer-Cartan equations, we have $g(\alpha_1+\alpha_2, \alpha_i)=g(F,e_i^-\wedge e_i^+)=g(F,J(e_i^+)\wedge e_i^+)=g(e_i^+,e_i^+)^2=G_i>0$. Now from $g(\alpha_i,J(\alpha_i))=0$ we obtain after a short calculation that $g(\alpha_1,\alpha_2)=-\frac{a}{b}g(\alpha_1,\alpha_1)=\frac{a}{c}g(\alpha_2,\alpha_2)$. Then from the previous equation we have
$$
\begin{array}{lllll}
g(\alpha_1+\alpha_2, \alpha_1)&=&\frac{b-a}{b}g(\alpha_1,\alpha_1)&=&G_1\\
g(\alpha_1+\alpha_2, \alpha_2)&=&\frac{a+c}{c}g(\alpha_2, \alpha_2)&=&G_2
\end{array}
$$
From the fact that $G_i$ and $g(\alpha_i, \alpha_i)$ are positive and $b,c\neq 0$ we obtain that $b(b-a)>0, \frac{a^2+1-ab}{a^2+1}>0$. This is necessary condition for the existence of homogeneous CYT structure on $SU(2)\times SU(2)$. Moreover one can check that it is also sufficient because we can determine the metric to be the form $g=g_{\alpha}+ \sum G_i(|e_i^+|^2+|e_i^-|^2)$ where $g_{\alpha}$ is Hermitian metric on $Span(\alpha_1,\alpha_2)$. This follows if we fix $g(\alpha_1,\alpha_1)=1$ and use the above equations to determine $G_i$ and $g_{\alpha}$.
The result can be formulated as
\begin{proposition}
The invariant complex structure $J$ on $SU(2)\times SU(2)$ determined by (\ref{cs}) admits a homogeneous CYT structure if and only if $b(b-a)>0, a^2+1-ab>0$.
\end{proposition}

We notice that for the space ${\it S}^1\times SU(2)$ any invariant complex structure carries an invariant CYT metric and the same is true for $T^3\times SU(2)$. This will be used in Section 5. It is a natural question to ask whether the complex structures on $SU(2)\times SU(2)$ which do not admit an invariant CYT metric admit any such metric. This is still an open question and the symmetrization above can not be applied.

\section{Compact complex parallelizable manifolds}
Another type of compact complex homogeneous manifolds with vanishing first Chern class
 are the complex parallelizable manifolds - i.e. the
manifolds with holomorphic parallelization of its holomorphic tangent bundle.

If the canonical bundle is holomorphically trivial, then from formula (\ref{chern}) follows that $\rho =
dd^C\log |\Omega|$ for a nonvanishing holomorphic section $\Omega$ and
\begin{equation}\label{paraleliz}
\rho^B = dd^c\log |\Omega| - d\delta F
\end{equation}

It is well known that complex paralelizable manifolds are of the form $G^c/\Gamma$, where $G^c$ is a complex Lie
group and $\Gamma$ is a cocompact lattice. For the Hermitian geometry of such manifolds there is the following
result which is easily deduced from the main theorem in  \cite{AG}:

\begin{theorem} For a compact complex parallelisable manifold any left invariant metric is a balanced metric, i.e $\delta F=0$.
\end{theorem}

A sketch of the proof is as follows: The Levi-Civita connection satisfies $g((D_XJ)Y,Z)=g(ad_{JZ}X,Y)$ from the standard
formulas. Then $\delta F(X)=tr(ad_{JX})$, which does not depend on $g$. But a Lie group $G$ admits a compact
quotient only if it is unimodular and for such groups $tr(ad_X)=0$ for any left-invariant $X$.

Now using it we can prove the following:
\begin{theorem} For a compact complex parallelisable manifold any left invariant metric is a CYT metric. Moreover all Ricci forms of the
canonical Hermitian connections vanish.
\end{theorem}

\bproof Let $M=G^c/\Gamma$ and let $\mathfrak{g}^c$ is the complex Lie algebra of $G^c$. Then there exist a holomorphic
left invariant vector fields $X_1, X_2, ... X_n$ on $G^c$ which form a basis for ${\bf g}^{(1,0)}$ and
 moreover
$[X_i,X_j]=c^k_{ij}X_k, [X_I,\overline{X_j}]=0$. If we define the
(1,0)-forms
 $\alpha_i$ to be dual of $X_i$ i.e. $\alpha_i(X_j)=\delta_{ij},
\alpha_i(\overline{X_J}) = 0$, then because of the last identity,
$\alpha_i$
 constitute a holomorphic 1-forms. Now consider any left
 invariant metric $g$ on $G$ (and hence on $M$). The form $\Omega=\alpha_1\wedge \alpha_2\wedge ... \alpha_n$ is
 a holomorphic (n,0)-form with constant norm.
From (\ref{paraleliz}) and the previous theorem we conclude that all Ricci forms for the canonical connections
vanish. \eproof

${\it Examples:}$ There is a classification of all complex Lie algebras in (complex) dimension three \cite{OV}.
These are the following:

i) The abelian algebra ${\mathbb C}^3$

ii) The complex Heisenberg algebra which is two-step nilpotent.

iii) The direct sum ${\it s}_2({\mathbb C})\oplus{\mathbb C}$ of the 2-dimensional solvable algebra ${\it
s}_2({\mathbb C})$ defined by $[X_1,X_2]=X_2$ for the basis $X_1,X_2$ and a center ${\mathbb C}$.

iv) Two irreducible solvable Lie algebras defined as follows:
$$ s_3({\mathbb C}) = span\{ X_1,X_2,X_3 | [X_1,X_2]=X_2,[X_1,X_3]=X_2+X_3 \}$$
and
$$ s_{3,\lambda}({\mathbb C}) = span\{ X_1,X_2,X_3 | [X_1,X_2]=X_2,[X_1,X_3]=\lambda X_3 \}$$

v) $sl(2,{\Bbb C})$.

It is known that the Lie groups of the complex Heisenberg algebra, the solvable Lie algebra $s_{3,-1}$ and $sl(2, {\Bbb C})$
admit cocompact latices, the others do not since they are not unimodular. In particular one has a classification
of the 3-dimensional compact complex parallelizable manifolds.

\begin{remark} In \cite{Na} are considered the deformations of some of the compact
complex parallelisable manifolds above. In particular it is shown that there exists a deformation of a
parallelizable manifolds which provides non-parallelizable ones. These are again complex manifolds of vanishing
first Chern class. Whether they admit a CYT metric is an open question.
\end{remark}


\section{Compact complex homogeneous manifolds of dimension 3 with vanishing first Chern class}

 We use in this section the result due to Tits about the classification of the compact complex homogeneous
 3-manifolds.
It is described as follows (Theorem 6.3 in Tits \cite{T}) :
 Any compact complex homogeneous space of dimension 3 is one of the following:

$i)$ a compact complex parallelizable manifold classified in Section 4.

$ii)$ a torus fibration of appropriate fiber dimension over ${\bf CP}^2$,${\bf CP}^1\times {\bf CP}^1$,
 or ${\bf CP}^1$.

$iii)$ a generalized flag manifold

$iv)$ reducible compact complex manifold

From this list $i)$ $ ii)$ and possibly $iv)$ are with vanishing first Chern class. This follows form Proposition 2 for $i)$
and is standard for $ii)$.  Moreover in $iv)$ we have products of a two-dimensional torus and a homogeneous compact
complex surface. One can check (again in Tits \cite{T}) that the homogeneous compact complex surfaces with
vanishing first Chern class are the following:

$v)$ torus

$vi)$ Hopf surface fibered over ${\bf CP}^1$ with a homogeneous complex structure.

\noindent Combining the results of this and the previous sections we obtain:

\begin{theorem} Any compact complex homogeneous manifold of complex dimension 3 with vanishing first Chern class
is biholomorphic to one of the above. Except for some invariant complex structures on $SU(2)\times SU(2)$ (or its finite factor), such space admits a homogeneous CYT structure. The set of invariant complex structures on $SU(2)\times SU(2)$ which admit a homogeneous CYT structure is open and is described in Proposition 3.
\end{theorem}

The proof is a case by case argument where only the case of the torus bundles in $ii)$ above is not evident. But all of them are
finitely covered by $SU(2)\times SU(2)$ or the product of a 2-tori with a Hopf surface. These cases were considered in Section 3.3.

\begin{remark} Some real homogeneous
spaces which admit left invariant complex structure are not complex homogeneous in this context. For example the
product of 5 dimensional real Heisenberg group with $S^1$ admits a compact quotient $M$ which is real
homogeneous and will be used in Section 8. There is also a left invariant complex structure which descends to
$M$. However in our terminology $M$ is only complex locally homogeneous.
\end{remark}

\section{ Complex homogeneous manifolds with invariant volumes}

In \cite{G} D.Guan proved the following:
\begin{theorem}
Every compact complex homogeneous space with an invariant volume form is a principal homogeneous complex torus
bundle over the product of a projective rational homogeneous space and a parallelizable manifold.
\end{theorem}

 There is more information about such bundles in \cite{G}. It is shown that the
bundle $\pi: M\rightarrow G/K\times D$, where $G$ is a compact Lie group and $H$ is closed subgroup, such
that $G/K$ is a rational homogeneous space and $D$ is a complex parallelizable
space, arises as a factor of a product of two principal complex torus
bundles of the same rank. One is $\pi_1: G/H\rightarrow G/K$, which is the Tits fibration for
$G/H$ and the other is $\pi_2:D_1\rightarrow D$ where $D_1$ is again compact
complex parallelizable and the fiber is a complex torus, which is in the center
of $D_1$. So $M=G/H\times_{T^n} D_1$ for some torus $T^n$ acting in anti-diagonal manner.

The main goal of this section is to prove the following:

\begin{theorem}\label{main}
Every compact complex homogeneous space with invariant  volume and a vanishing first Chern class is a principal
torus bundle over a product of a standard homogeneous CR manifold and a complex parallelizable manifold. After a
homogeneous complex deformation it admits a homogeneous CYT structure.
\end{theorem}
The standard homogeneous CR space here is obtained from the principle circle bundle over a generalized flag manifold with
 characteristic class equal to the first Chern class of the flag space, divided by its multiplicity. Such fibre bundle is called Boothby-Wang fibration.
 By homogeneous CYT structure we mean a Hermitian metric invariant under the action of some transitive Lie group of
biholomorphisms, and a homogeneous complex deformation is a 1-parameter family of homogeneous complex structures.
Before we start the proof of the theorem we need the following two lemmas:
\begin{lemma}
For any $(a_1,a_2,...,a_k)\in {\Bbb Z}^k$ with $gcd(a_1,...,a_k)=1$, there are
$(k-1)$ elements of $\Bbb{Z}^k$ which together with $(a_1,a_2,...,a_k)$
determine a matrix in $SL(k, \Bbb{Z})$
\end{lemma}

We sketch the proof of the Lemma for sake of completeness. We use the Bezout
Lemma, that for any $(a_1,a_2)$ there
 are $(x,y)$ with $gcd(x,y)=1$, such that $a_1x-a_2y =gcd(a_1,a_2)$. Start
 with the standard basis $(e_1,...,e_k)$ of $\Bbb{Z}$ and change it first to
 $e_1'=(\frac{a_1}{gcd(a_1,a_2)}e_1+\frac{a_2}{gcd(a_1,a_2)}e_2,
 e_2'=xe_1+ye_2,e_3'=e_3,...,e_k'=e_k)$. Then similarly we change the second
 basis to a new one using the solution of the equation
 $gcd(a_1,a_2)x_2-a_3y_2=gcd(a_1,a_2,a_3)$ with
 $e_1"=\frac{gcd(a_1,a_2)}{gcd(a_1,a_2,a_3)}e_1'+\frac{a_3}{gcd(a_1,a_2,a_3)}e_3'$
$e_2"=e_2'$ and $e_3"=x_2e_1'+y_2e_3'$. We see that at each step the matrix of
the basis is in $SL(k, \Bbb{Z})$ and we can continue until the $k$th step will
produce the necessary basis.

\begin{lemma}\label{SL}
Let $M\rightarrow B$ be a principal complex torus bundle with real
characteristic classes $\{\omega_1, \omega_2,...\omega_{2k}\}$ over a complex
base $B$, such that $M$ carries a torus invariant complex structure. If $A =
(a_{ij})$ is a matrix in $SL(2k,\Bbb{Z})$, then the classes $\omega_i' =\sum_j
a_{ij}\omega_j$ determine an equivalent torus bundle equipped with an equivalent
complex structure. In particular one can choose the classes in such a
way that the first $m$ are ${\Bbb Z}$-independent and the rest $2k-m$ are zero.
\end{lemma}
\bproof
 Let the torus action
 which induces the toric bundle structure on $M$ with characteristic classes
$(\omega_1,...,\omega_{2k})$ be given by $(e^{i\theta_1}, ...,
e^{i\theta_{2k}}):M \rightarrow M$. There are connection forms
$\eta_1,...,\eta_{2k}$, such that the curvatures are forms in
$(\omega_1,...,\omega_{2k})$ and the complex structure on $M$ is given by the
standard construction (see \cite{GGP}). Then the torus action determined
by $(e^{i\theta'_1},...e^{i\theta'_{2k}})$, where
$\theta'_i=\sum_ja_{ij}\theta_j$ is free. Topologically the two actions
determine the same fibers and the same base. Now consider the forms
$\eta_1',...,\eta_{2k}'$ given by $\eta'_i=\sum_ja_{ij}\eta_j$ which provide a
$T^{2k}$-connection in the second bundle. Then the characteristic classes are
$(\omega_1',...,\omega_{2k}')$ and the horizontal spaces are the same. It also
induces an automorphism of the complex structure on the fiber, since the matrix
$(a_{ij})$ is in $SL(2k,\Bbb{Z})$ and corresponds to an automorphism of the
complex torus (the fiber) induced by an automorphism of its underlying lattice.
So the two complex structures on $M$ are equivalent. For the second claim we
apply the previous Lemma. \eproof

{\it Proof of Theorem \ref{main}:} From D.Guan's result any compact complex homogeneous space with invariant volume
admits a complex torus fiber bundle structure as $M\rightarrow G/K\times D$ with characteristic classes
$(\omega_1+\alpha_1, \omega_2+\alpha_2,...\omega_{2k}+\alpha_{2k})$, where $(\omega_1,...,\omega_{2k})$ are the
characteristic classes of $G/H\rightarrow G/K$ which are of type (1,1) and $(\alpha_1,...,\alpha_{2k})$ are the
characteristic classes of $D_1\rightarrow D$.  Notice that by the properties of the averaging map $\mu$ in Section 4 there is a unique $G$-invariant representative in each class $\omega_k$. The second fibration is a complex torus fibration and its
(complex) characteristic classes are of type (2,0) with respect to the complex structure on the base $D$. In
particular $\alpha_i$ are of type (2,0)+(0,2), i.e they have representatives of this type. From now on we denote
the invariant representatives with the same letters. Now the first Chern class of $M$ vanishes iff
$c_1(G/K\times D)$ is represented as a ${\Bbb Z}$-combination of the characteristic classes of $M\rightarrow
G/K\times D$. Since $c_1(D)= 0$,  there are integer constants $(d_1, d_2,...,d_{2k})$ such
that $d_1\alpha_1+...+d_{2k}\alpha_{2k} = 0$ and $d_1\omega_1+...+c_{2k}\omega_{2k}=\rho$, where $\rho=c_1(G/K)$
is the class of the Koszul form of $G/K$. This follows from the type decomposition of $(\alpha_1, ..., \alpha_{2k}, \omega_1, ..., \omega_{2k})$. Moreover since $\alpha_1,...,\alpha_{2k}$ are real and imaginery parts
of the (2,0)-forms $\alpha_{2i+1}+i\alpha_{2i+2}$, then $\alpha_{2i+2}(JX,Y)=\alpha_{2i+1}(X,Y)$ for any tangent
vectors $X,Y$. From here we have $-d_2\alpha_1+d_1\alpha_2-d_4\alpha_4+...+d_{2k-1}\alpha_{2k} = 0$, Now we
change the basis for the characteristic classes of the fibration $M\rightarrow G/K\times D$ by choosing
$m\omega_1'=d_1(\omega_1+\alpha_1)...+d_{2k}(\omega_{2k}+\alpha_{2k})= \rho$, where $m = gcd(d_1,...d_{2k})$.
Then we complete $\omega_1'$ to a basis of characteristic classes $(\omega_2'+\alpha_2',
...,\omega_k'+\alpha_k')$ which determine the same bundle $M$ such that $\omega_i'$ are on $G/K$ and $\alpha_i'$
are on $D$.
 Since $\omega_1'$ is proportional to the first Chern class
of $G/K$, it defines a bundle $G/H_1\rightarrow G/K$ which is a standard
homogeneous CR space and the bundle structure is the Boothby-Wang
fibration. To obtain a CYT metric we use the formula (\ref{torusrel}) for the
relation of the Ricci forms in generic torus bundle $X\rightarrow B$ over Hermitian
manifold $B$:
$$
\rho^B_X = \pi^*(\rho^B_B - \sum g(F,\omega_i)\omega_i)
$$
when $g(F,\omega_i)=constant$ and $F$ is the fundamental form on $B$. Since $g(F,\alpha_i)=0$ we need to find
$g$ on $G/K\times D$, such that $\rho^B_{G/K\times D}= \sum g(F,\omega_i)(\omega_i+\alpha_i)$. Because
$\alpha_i$ is of type (2,0)+(0,2), we have to find $g$ and a basis of characteristic classes, such that
$g(F,\omega_i)=0$ for $i>1$. Now consider the bundle $G/H\rightarrow G/K$ from the Guan's construction with
characteristic classes $(\omega_1,...,\omega_{2k})$. Then $G/H$ is itself a complex manifold of vanishing first
Chern class and the same bundle is determined also by some classes $(\omega_1',...,\omega_{2k}')$, with
$\omega_1'=(1/m) c_1(G/K)$ by Lemma \ref{SL}.
Then by Theorem \ref{compact} the complex structure on $G/H$ can be deformed to one that admits a homogeneous CYT metric $h$. The deformation
preserves the Tits fibration and doesn't change the complex structure on $G/K$. In particular $\rho^B_{G/K}=\sum
g(F,\omega_i')\omega_i'$ for the fundamental form $F$ of the metric $g$ induced on $G/K$ from $h$. However we
can choose the classes $\omega_i'$ to be such that part of them are zero, and the rest are linearly independent
and represented by invariant forms. Since $\rho^B$ is proportional to $\omega_1'$, all $g(F,\omega_i')=0$
for $i\geq2$ and $\rho^B=g(F,\omega_1')\omega_1'$. Then the metric we are looking for on $G/K\times D$ is a
product of $g$ and any invariant Hermitian metric on $D$. The above mentioned toric bundle construction produces
a CYT metric on $M$.

The last step of the proof is to show that the metric could be chosen so  that it is homogeneous. We noticed that $\omega_i$ are chosen $G$-invariant. For the second factor $D$, there is a Lie group $H$ acting transitively, such that the isotropy subgroup is a cocompact lattice. Then by the averaging which is used in \cite{GF} one can find $H$-invariant representatives for $\alpha_i$. This will produce a metric which is invariant
under the action of some real transitive Lie group obtained as an abelian extension of $G\times H$ which will also leave the volume form invariant.

\eproof

\section{Topological properties}

In this section we discus the integral cohomology of the space  $SU(4)/U(1)$ from Proposition \ref{example1}, although our
considerations work for $SU(2n)/U(1)$ in general. The calculations follow \cite{Esch} and are based on the classical methods developed by Borel and Serre, see  \cite{Borel}. Let $U(1)$ be the subgroup defined by  $(e^{2k\pi it}, e^{2l\pi it}, e^{2m\pi it}, e^{-2(k+l+m)\pi it})$
 in $SU(4)$. Assume that $gdc(k,l,m)=1$ so that the action is free. The first Chern class vanishes if $m=-3k-2l$, or $U(1)$ acts with weights $(k,l,-3k-2l,2k+l)$ by Proposition \ref{example1}. The embedding induces a projection of the corresponding classifying
spaces $\rho:BU(1)\rightarrow BSU(4)$. The integral cohomology $H^*(BU(1))$ are identified
 with the ring of polynomials $\mathbb{Z}[s]$ of one variable of degree 2. Also $H^*(BSU(4))$
is a subring of $H^*(BU(1)^4)$ generated by the elementary symmetric polynomials in 4 variables with the
restriction that the first one vanishes. More precisely, if $(t_1,t_2,t_3,t_4)$ are the independent variables
corresponding to the generators of $H^*(BU(1)^4)$, then on $H^*(BSU(4))$, $t_4=-t_1-t_2-t_3$ and the generators
are $$u_4= -t_1^2-t_2^2-t_3^2+t_1t_2+t_2t_3+t_1t_3  $$
$$u_6=-2t_1t_2t_3+t_1^2t_2+t_1t_2^2+t_1^2t_3+t_1t_3^2+t_2^q2t_3+t_2t_3^2$$
       $$ u_8=-t_1t_2t_3(t_1+t_2+t_3)$$

      As a result the induced map $\rho^*$ on the cohomology has the following characterization:
\begin{lemma}
The map $\rho^*: H^*(BSU(4))\rightarrow H^*(BU(1))$ is determined by $\rho^*(u_4)=N_4s^2, \rho^*(u_6)=M_6s^3,
\rho^*(u_8)=K_8s^4$ where $M_4=-k^2-l^2-m^2+kl+lm+kl, N_6=-2klm+kl(k+l)+lm(l+m)+km(k+m),K_8= -klm(k+l+m)$
\end{lemma}

Let $EU(1)\rightarrow BU(1)$ be the universal $U(1)$ bundle and $ESU(4)\rightarrow BSU(4)$ the universal
$SU(4)$-bundle. Consider the following diagram:

$$\begin{array}{ccc}
 ESU(4)\times SU(4)& \rightarrow &SU(4)\\
\downarrow& &\downarrow\\
M&\rightarrow &X
\end{array}
$$

Then right column is the principal bundle defining our homogeneous space $X=SU(4)/U(1)$. The low horizontal
arrow is a projection with fibre the aciclic space $EU(1)\cong ESU(4)$ so determines an isomorphism on cohomology.
Then the space $M=(EU(1)\times SU(4))/U(1) = (ESU(4)\times SU(4))/U(1)$ is included in the diagram:

$$\begin{array}{ccc}
 ESU(4)\times SU(4)& \rightarrow & EU(1)\\
\downarrow& &\downarrow\\
M&\rightarrow &BU(1)
\end{array}$$

The fibres of the horizontal arrows could be identified with $SU(4)$. According to Borel, the
spectral sequence $E$ of $M\rightarrow BU(1)$ has second term $E_2=H^*(BU(1))\otimes H^*(SU(4))$. The ring
$H^*(SU(4))$ is $\Lambda(u_3,u_5,u_7)$ such that $(u_4,u_6,u_8)$ are images under the transgression of $(u_3,u_5,u_7)$ in the universal bundle
$ESU(4)\rightarrow BSU(4)$.

\begin{lemma}\label{spect}
The differentials of the spectral sequence $E$ of the bundle $M\rightarrow BU(1)$ are determined by $$d_{j}(1\otimes u_i)=0, j\leq i$$
$$d_{j+1}(1\otimes u_j)=\rho^*(u_{j+1})\otimes 1$$
\end{lemma}
\bproof
Consider the diagram:
$$\begin{array}{ccc}
M=ESU(4)\times SU(4)/U(1)& \rightarrow & ESU(4)\\
\downarrow& &\downarrow\\
BU(1)&\rightarrow &BSU(4)
\end{array}$$
The projection of the lower row is $\rho$ and it is covered by a map which induces the identity on the cohomology of the fibers.
The nontrivial differentials of the universal bundle $ESU(4)\rightarrow BSU(4)$ are given by: $$d_{j}(1\otimes u_i)=0, j\leq i$$
$$d_{j+1}(1\otimes u_j)=u_{j+1}\otimes 1$$
Then the Lemma follows by naturality. \eproof

From the previous two lemmas one can use $E$ to compute the cohomology of $M$ and
hence of $X$. First notice that $E_2=E_3=E_4$. Since $d_4(1\otimes u_3)=M_4s^2\otimes 1$ by Lemma \ref{spect} it follows that  $Im
d_4=<M_4s^2>\otimes 1$ and $Ker d_4 = H^*(BU(1))\otimes<1,u_5,u_7>$ where the last notation means the ideal generated by the
decomposable elements with first factor in $H^*(BU(1))=\mathbb{Z}[s]$ and second factor one of the $1, u_5$ or $u_7$. Hence $E_5 = Ker d_4/Im d_4 = \mathbb{Z}[s]/(M_4s^2)\otimes<1,u_5,u_7>$. Since again
$E_5=E_6$, the next step to consider is $d_6(1\otimes u_5)=N_6s^3\otimes1$. Now assume that
$gcd(M_4,N_6)=L$. Then $1\otimes (M_4/L)u_5 \in Ker d_6$ and $d_6$ acts on $E_6$. So $Ker d_6 =
\mathbb{Z}[s]/(M_4s^2)\otimes <1,(M_4/L)u_5,u_7>$ and $Im d_6=<N_6s^3\otimes 1>$. By the assumption, $<M_4s^2,
N_6s^3>=<M_4s^2,Ls^3>$, since there are $a$ and $b$ such that $aM_4+bN_6=L$. Then we have $Ker d_6/Im d_6 =
\mathbb{Z}[s]/<M_4s^2,Ls^3>\otimes<1,(M_4/L)u_5,u_7>=E_7=E_8$. We can also see that $d_8(1\otimes u_7)=K_8s^4\otimes 1 =
0$ in $E_8$ if $gcd(L,K_8)=gcd(M_4,N_6,K_8)=1$. Then all other differentials vanish and $E_{\infty}=E_7$. Now if we denote by $w=s\otimes 1,
v_5=1\otimes (M_4/L)u_5, v_7=1\otimes u_7$, then we have:

\begin{theorem}\label{cohom}
Let $k,l$ are relatively prime such that $M_4=13k^2+7l^2+16kl, N_6=6k^3+2l^3+26k^2l+18kl^2, K_8= 6k^3l+2kl^3+7k^2l^2$ are also relatively prime and $X=SU(4)/U(1)$ is given by the action of $U(1)$ with weights $(k,l,-3k-2l,2k+l)$. The cohomology ring $H^*(M, \mathbb{Z})$ is generated by $w\in H^2, v_5\in H^5, v_7\in H^7$ with the relations
$(M_4/L)w^2=w^3=w^2v_5=w^2v_7=0$ where $L=gcd(M_4,N_6)$. In particular the nonzero cohomology are $H^2=\mathbb{Z}, H^4=\mathbb{Z}_{|M_4/L|},
H^7=\mathbb{Z}^2, H^0=H^5=H^9=H^{12}=H^{14}=\mathbb{Z}$
\end{theorem}

The choice of the weights is such that $c_1(X)=0$. From here and the condition $c_1=0$, we see that all rational Chern numbers and hence all rational Pontriagin numbers of $M$ vanish. A simple considerations regarding the relations of the characteristic classes in the Tits fibration show that only $c_2$ and $p_1$ could be nonzero but torsional. When we vary $k$ and $l$ we obtain:
\begin{corollary}
There are infinitely many simply-connected non-flat CYT manifolds of dimension 14 with the same Hodge and Betti numbers and torsional Chern classes.
\end{corollary}
\bproof
The curvature of the space $SU(4)/U(1)$ is the curvature of the canonical connection which is given in \cite{KN}. It vanishes on $SU(4)$ but one can check directly that is nonzero on the quotient. The statement for the Betti numbers follows from Theorem \ref{cohom}. The Hodge numbers for such quotients are calculated in \cite{Gr}, Proposition 5.1 and one can check directly that they are equal. The only part which remains to be proven is that there are infinitely many choices of numbers $k$ and $l$ as in the Theorem \ref{cohom} for which $|M_4/L|$ is unbounded. First take $l=1$ and notice that $S=gcd(M_4,K_8)=gcd(13k^2+16k+7,6k^3+7k^2+2k)$ divides $5k+16$. There are infinitely many $k$ for which this number is prime, which easily gives that $S=1$. Then $gcd(M_4, N_6, K_8)=1$ and the conditions of Theorem \ref{cohom} are met. Similarly one can check that $L$ divides a linear polynomial of $k$, so $M_4/L$ will grow infinitely.
\eproof

\begin{remark} The Dolbeaut cohomology of the compact complex parallelizable manifolds are  subject of
intensive investigation. In general $H^{(p,q)}(G/\Gamma)=\Lambda^q(g^+)\otimes H^{(0,q)}(G/\Gamma)$. where $g^+$
is the $+\sqrt{-1}$-eigenspace of the Lie algebra of $G$. The second factor depends both on $G$ and $\Gamma$. In
case $G=SL(2,\mathbb{C})$ there are different $\Gamma$ which lead to quotients with different cohomology rings
and different types of Kuranishi spaces of deformations \cite{G}.
\end{remark}

\section{ Relation to the Strominger's equations in heterotic string theory}

 In 1986
A.Strominger \cite{Strom} analyzed the heterotic superstring backgrounds with spacetime supersymmetry. His model is
based on Hermitian manifolds which are CYT spaces with holomorphically trivial canonical bundle. In terms of the
Hermitian geometry it consists of conformally balanced complex 3-manifold with holomorphic (3,0)-form of
constant norm and an anomaly cancelation condition. The manifold is endowed with an auxiliary semistable bundle
with Hermitian-Einstein connection $A$ with curvature $F_A$ and the anomaly cancelation condition is:
$$
dH = 2i\partial\overline{\partial}F= dd^c F = \frac{\alpha'}{4}[tr(R\wedge R)-
tr(F_A\wedge F_A)]
$$
 for the K\"ahler form $F$. Here $R$ is the curvature of some metric connection.  In fact {\it any} metric connection with curvature $R$ solving the anomaly cancelation condition with $\alpha'>0$ leads to a physically meaningful solution as noted in \cite{BBFTY}. The first solutions on non-K\"ahler manifolds of this system were constructed
only recently by J.Fu and S.T.Yau \cite{yau}. \vspace{.1in}

We start here with solutions of the Strominger's system in which
the anomaly cancellation equation have trivial instanton $F_A=0$. Later we provide also a nontrivial instanton solutions.
Let $\{e^1,Je^1,e^2,Je^2,e^3,Je^3\}$ be an ordered  unitary co-basis for a
complex structure $J$ and Hermitian metric $g$ on a vector space. Consider the nilpotent
Lie algebra whose sole non-trivial structure equation is as follows:
\[
d(Je^3)=e^1\wedge Je^1-e^2\wedge Je^2.
\]
With the given ordered basis, we could consider this algebra as
$(0,0,0,0,0,12-34)$. If we change the order of the basis to
$\{e^1,Je^1,Je^2,e^2,e^3,Je^3\}$, we obtain the "standard basis" in
classification of six-dimensional nilpotent algebras. The structure equation is
$(0,0,0,0,0,12+34)$. As a real Lie algebra, it is the direct sum of a real
five-dimensional Heisenberg algebra and a one-dimensional trivial algebra.

Since $d(e^3+iJe^3)\in \Lambda^{(1,1)}$ and the exterior differential of
$e^1+iJe^1$ and $e^2+iJe^2$ are equal to zero,  the complex structure is
integrable. As $d\Lambda^{(1,0)}\subset \Lambda^{(1,1)}$, the complex structure
is also \it abelian. \rm

A direct calculation shows that:
\begin{eqnarray*}
F^2&=&-2\sum_{1\leq j<k\leq 3} \left((e^j+iJe^j)\wedge
(e^j-iJe^j)\wedge(e^k+iJe^k)\wedge (e^k-iJe^k)\right)\nonumber\\
&=&8\sum_{1\leq j<k\leq 3}\left( e^j\wedge Je^j\wedge e^k\wedge Je^k\right).
\end{eqnarray*}
Therefore,
\begin{eqnarray*}
dF^2&=& 8\sum_{1\leq j<k\leq 3} d\left( e^j\wedge Je^j\wedge
e^k\wedge Je^k\right)\nonumber\\
&=& -8 \left( e^1\wedge Je^1\wedge e^3\wedge dJe^3+e^2\wedge
Je^2\wedge e^3\wedge dJe^3\right)\nonumber\\
&=&-8 \left( -e^1\wedge Je^1\wedge e^3\wedge(e^2\wedge Je^2)
+e^2\wedge Je^2\wedge e^3\wedge (e^1\wedge Je^1)\right)\nonumber\\
&=&0
\end{eqnarray*}
so the metric is balanced. Moreover the form
$(e^1+iJe^1)\wedge(e^2+iJe^2)\wedge(e^3+iJe^3)$ is easily seen to be closed and
holomorphic.

Next note that since $d(Je^3)$ is a type (1,1)-form, we have $Jd(Je^3)=dJe^3$.
Therefore,
\begin{eqnarray*}
 dd^cF &=& dJdJF =dJdF\\
&=&2 dJd(e^3\wedge Je^3)=2 dJ(-e^3\wedge dJe^3)\\
&=&-2d(Je^3\wedge JdJe^3)\\
&=&-2d(Je^3\wedge dJe^3)\\
&=&-2(dJe^3)\wedge (dJe^3).
\end{eqnarray*}
Now we take a connection form $\omega$ to be
\[
\omega = \left(\begin{array}{llllll}
0&aJe^3&0&0&0&0\\
-aJe^3&0&0&0&0&0\\
0&0&0&0&0&0\\
0&0&0&0&0&0\\
0&0&0&0&0&0\\
0&0&0&0&0&0
\end{array}\right).
\]
with $a=constant$. Then we obtain $tr(R\wedge R)= -a^2(dJe^3)^2 = +a^2dd^cF =
2a^2i\partial\overline{\partial} F$. So this provides a solution to the
Strominger's system.

 Now we can also construct an instanton bundle which still satisfies the
 anomaly cancellation condition. We take any Hermitian Yang-Mills connection on
 a holomorphic vector bundle on the base 4-tori. To ensure its existence we may
 assume that the base is an abelian surface. Then we can take its pull-back to
 $M$. The curvature satisfies $tr(F\wedge F) = b e^1\wedge Je^1\wedge e^2\wedge
 Je^2$ for some function $b$. If we take the bundle to be homogeneous,
 $b=constant$. So $dd^cF= (\frac{1}{a^2}-b)e^1\wedge Je^1\wedge e^2\wedge
 Je^2$. By choosing $a$ sufficiently small, we have a solution for the anomaly
 cancelation condition with positive $\alpha'$.

\begin{remark}
The solution of \cite{yau} uses the Chern connection and produces a nonconstant dilaton. The connection is
chosen because the curvature term has appropriate type (2,2). In physics
a preferred choice of connection is a metric non-Hermitian connection with
skew-symmetric torsion equal to the negative of the torsion of the Bismut
connection. Unfortunately the example above doesn't provide a solution for this
connection since for its curvature $R$ one obtains $tr(R\wedge R)=0$.

The compact factors of $H_5\times \mathbb{R}$ are toric bundles over 4-tori. For such
there is a vanishing theorem in  \cite{BBFTY}
whenever the contribution of $tr(R\wedge R)$ is zero.
After the results of this section were reported the paper \cite{FIM} appeared, where similar solutions of the Strominger's equations were given on many 6-dimensional nilmanifolds, using different invariant connection and instanton. Their solutions, like our have constant dilaton and instanton with vanishing Euler number. From physics point of view more realistic models require this Euler number to be six.
Such non-K\"ahler solutions are still unknown. It is worth to notice that the 3-dimensional compact complex parallelizable manifolds admit solutions of the Strominger's equations, but with $\alpha'<0$. The role of such solutions is unclear and we don't provide any details here.
\end{remark}

{\bf Acknowledgments:} This paper started as a joint project with Y.S. Poon and the author is particularly
grateful to him for the numerous discussions and help. Special thanks go to D. Alekseevsky, who  pointed out the
importance of the Koszul form at the initial stage of the project. I would like to thank the referees for their remarks which improved the presentation in many places.

\noindent Department of Mathematics

\noindent Florida International University

\noindent Miami, FL 33199

\noindent grantchg@fiu.edu.

\end{document}